\input amstex
\documentstyle{amsppt}
\magnification 1200
\NoRunningHeads
\NoBlackBoxes
\document

\def\v{{\bold v}}
\def\w{{\bold w}}
\def\G{{\Bbb G}}

\def\ell{{\text{ell}}}

\def\qdet{\text{qdet}}
\def\k{\kappa}
\def\RR{\Bbb R}
\def\be{\bold e}
\def\bR{\overline{R}}
\def\tR{\tilde{\Cal R}}

\def\R{\Cal R}
\def\h1{\hat{\bold 1}}
\def\hV{\hat V}

\def\hz{\hat \z}
\def\hV{\hat V}

\def\Ind{\text{Ind}}

\def\Hom{\text{Hom}}

\def\Ua{U_q(\tilde\g)}
\def\U2{{\Ua}_2}
\def\g{\frak g}

\def\d{\partial}

\def\z{\bold z}
\def\Id{\text{Id}}
\def\<{\langle}
\def\>{\rangle}
\def\o{\otimes}
\def\e{\varepsilon}

\def\Id{\text{Id}}
\def\End{\text{End}}

\topmatter
\title Quantization of Lie bialgebras, IV
\endtitle
\author {\rm {\bf Pavel Etingof and David Kazhdan} \linebreak
\vskip .1in
Department of Mathematics\linebreak
Harvard University\linebreak 
Cambridge, MA 02138, USA\linebreak
e-mail: etingof\@math.harvard.edu\linebreak kazhdan\@math.harvard.edu}
\endauthor
\endtopmatter

\head Introduction \endhead

This paper is a continuation of \cite{EK3}. In \cite{EK3}, we introduced 
the Hopf algebra $F(R)_\z$ associated to a quantum R-matrix 
$R(z)$ with a spectral parameter defined on a 1-dimensional 
connected algebraic group $\Sigma$, and a set of points
$\z=(z_1,...,z_n)\in \Sigma^n$. This algebra is generated by 
entries of a matrix power series $T_i(u)$, $i=1,...,n$, subject to 
Faddeev-Reshetikhin-Takhtajan type commutation relations, and 
is a quantization of the group $GL_N[[t]]$. 
 
In this paper we consider the quotient $F_0(R)_\z$ 
of $F(R)_\z$ by the relations $\qdet_R(T_i)=1$, where 
$\qdet_R$ is the quantum determinant associated to $R$
(for rational, trigonometric, or elliptic R-matrices). 
This is also a Hopf algebra, which is a quantization of the group 
$SL_N[[t]]$.  

This paper was inspired by \cite{FR}. The main goal of this paper 
is to study the representation theory of the algebra $F_0(R)_\z$ and
of its quantum double, and show how the consideration of
coinvariants of this double (quantum conformal blocks) 
naturally leads to the quantum 
Knizhnik-Zamolodchikov equations of Frenkel and Reshetikhin \cite{FR}. 
Our construction for the rational R-matrix 
is a quantum analogue of the standard derivation of the Knizhnik-Zamolodchikov 
equations in the Wess-Zumino-Witten model of conformal field theory 
\cite{TUY}, 
and for the elliptic R-matrix is a quantum analogue of the construction of 
\cite{KT}. 

Our result is a generalization of the construction of Enriques 
and Felder \cite{EF}, which appeared while this paper was in preparation.
Enriques and Felder gave a derivation of the quantum KZ equations
from coinvariants in the case of the rational R-matrix and N=2. 

The results of this paper for the rational R-matrix (the Yangian case) 
can be directly generalized to the case of any simple Lie algebra $\g$ 
(what we do here corresponds to $\g={\frak{sl}_N}$). We did not include 
this generalization here since for a general $\g$ it is more difficult 
to write explicit formulas. 

We note that this paper does not use the results from \cite{EK1,EK2} 
on the existence of quantization. 

Finally, we would like to explain the relationship between 
the present paper and the papers \cite{FR,KS}, which are devoted to the same 
subject. The papers 
\cite{FR,KS} generalize to the quantum case the construction of 
Tsuchiya-Kanie (\cite{TK}), which represents conformal blocks as intertwiners
between a highest weight representation and the tensor
 product of a highest weight 
and a finite dimensional representation of an affine Lie algebra.
This allows to obtain the quantum KZ equations, but does not allow 
to define the quantum analogue of fusion of highest weight modules  
for affine Lie algebras. On the contrary, the approach of the present paper
generalizes the coinvariant approach \cite{TUY}, and allows to define
a quantum analogue of fusion (i.e. of vertex operator calculus). 

However, the results of our paper hold only for formal quantum parameter $h$, 
while the results of \cite{FR,KS} are valid for numerical $h$. 
It is an interesting question what is the analogue of our results for 
numerical $h$.  

In the next paper we will discuss the notion of a quantum vertex operator 
algebra, which naturally arises from the setting of 
``Quantization of Lie bialgebras III and IV'', and describe the 
quantum analogue of the affine vertex operator algebra. 

\head Acknowledgements\endhead

The authors were supported by NSF grant DMS-9700477. 
  
\head 1. The algebra $F_0(R)$\endhead 

\subhead 1.1. R-matrices with crossing symmetry \endsubhead

We recall the setting of \cite{EK3}. 

Let $k$ be an algebraically closed field of characteristic $0$. 

Recall \cite{EK3} that $R_{rat}(u),R_{tr}(u),R_{ell}(u)$ 
denote the rational, 
the trigonometric, and the elliptic R-matrix respectively.  
For example, $R_{rat}(u)=1-\frac{h(\sigma-1/N)}{N(u-h/N^2)}$,
where $\sigma$ is the permutation of components.

Let $R(u)\in \End(k^N\o k^N)((u))[[h]]$ be 
$R_{rat}(u/N),R_{tr}(u/N)$ or $R_{ell}(u/N)$. 

For any connected 1-dimensional algebraic group $\Sigma$ over $k$
(i.e. $\G_a,\G_m$, or an elliptic curve), we let 
$u$ be a canonical formal parameter on $\Sigma$ near the origin
(it is defined uniquely up to scaling). Given $R$, let 
$(\Sigma,u)$ be such that $R\in \End(k^N\o k^N)\o k(\Sigma)[[h]]$ 
and has a trivial stabilizer in $\Sigma$. 

As we mentioned in \cite{EK3}, the function $R(u)$ satisfies the 
quantum Yang-Baxter equation
$$
R^{12}(u_1-u_2)R^{13}(u_1-u_3)R^{23}(u_2-u_3)=
R^{23}(u_2-u_3)R^{13}(u_1-u_3)R^{12}(u_1-u_2),\tag 1.1
$$
and the condition
$$
\lim_{s\to 0}R(su,sh)=R_{rat}(u/N)=
1-\frac{h(\sigma-1/N)}{u-h/N}.\tag 1.2
$$ 

Besides, it is known that the function $R(u)$ 
has the property of crossing symmetry (see \cite{Ha} for the elliptic case):
$$
(((R(u)^{-1})^{t_1})^{-1})^{t_1}=
(((R(u)^{-1})^{t_2})^{-1})^{t_2}=
g(u)R(u+Nh).\tag 1.3
$$
where $t_i$ denotes transposition in the i-th component, and 
$g\in 1+hk((u))[[h]]$ 
is a scalar function. We will explain the meaning of this condition 
in section 1.4.

\subhead 1.2. The Hopf algebra $F(R)$\endsubhead

Recall from \cite{EK3} the definition and properties 
of the algebra $F(R)$. 

\proclaim{Definition} The algebra $F(R)$ is 
the h-adic completion of the 
algebra over $k[[h]]$ whose generators are the entries
of the coefficients of formal series $T(u)^{\pm 1}\in \End(k^N)\o F(R)[[u]]$,
 $T(u)=T_0+T_1u+...$, 
and the defining relations are
$$
\gather
T(u)T(u)^{-1}=T(u)^{-1}T(u)=1,\\
R^{12}(u-v)T^{13}(u)T^{23}(v)=T^{23}(v)T^{13}(u)R^{12}(u-v),
\tag 1.4\endgather
$$
\endproclaim

As we saw in \cite{EK3}, 
this algebra is a (topological) Hopf algebra, with the coproduct, 
counit, and antipode defined by
$$
\Delta(T(u))=T^{12}(u)T^{13}(u),\ \e(T(u))=
1,\ S(T(u))=T^{-1}(u).\tag 1.5
$$

Let $r$ be defined by $R=1-hr+O(h^2)$, and $GL_N(r)$
be the corresponding (infinite-dimensional) Poisson group 
defined in Chapter 2 of \cite{EK3}. 
By Proposition 
3.9(b) of \cite{EK3}, equations (1.1) and (1.2) imply 
that the Hopf algebra $F(R)$ 
is a quantization of $GL_N(r)$. 

\subhead 1.3. The quantum determinant for $F(R)$\endsubhead

In \cite{EK3}, we introduced the quantum determinant for the dual Yangian.  
Let us generalize this construction to the trigonometric and elliptic cases. 

Recall that for a coalgebra $A$, a right comodule over $A$ is 
a vector space $V$ with a linear map $\pi^*:V\to A\o V$ such that 
$(1\o \pi^*)\circ \pi^*=(\Delta^{op}\o 1)\circ \pi^*$. 

Let $E(u)$ be the basic comodule for $F(R)$ over $k[[u]]$. 
Namely, let $E(u)$ be the space $k^N[[u]]$ with the coaction 
$\pi^*_{E(u)}(\v )=T^{21}(u)(1\o \v )$. 
(it is easy to check that this formula indeed defines 
a right comodule). 

Consider the tensor product of comodules 
$E(u-\frac{h(N-1)}{2})\o 
E(u-\frac{h(N-3)}{2})..\o E(\frac{u+h(N-1)}{2})$. It follows from 
the representation theory of $F(R)$
 that this comodule has 
a unique 1-dimensional subcomodule $QDet$
(it is enough to use the functor from comodules to modules
defined in Section 2.3, and then apply the results of \cite{Ch}for the 
elliptic case). 
Denote the generator of $QDet$ by 
$\v_0$. Then we have $\pi^*_{QDet}(\v_0)=q(u)\o \v_0$, where 
$q(u)\in F(R)[[u]]$. The element $q(u)$ is central and grouplike, and it 
is called the quantum determinant. We will denote it by $\qdet_R(T(u))$. 
For the case of the rational R-matrix, the quantum determinant 
was discussed in Section 3.3 of \cite{EK3}. 

Let us write down an explicit expression for the quantum determinant. 
Let $\be_i$ be the standard basis of $k^N$. Let 
$\v_0=\sum C_{i_1...i_N}\be_{i_1}\o...\o \be_{i_N}$, and normalize 
$\v_0$ in such a way that $C_{1...N}=1$. Let 
$T=\sum E_{ij}\o T_{ij}$. Then 
$$
\qdet_R(T(u))=\sum_{i_1,...,i_N}C_{i_1,...,i_N}T_{1i_1}(u-h(N-1)/2)...
T_{Ni_N}(u+h(N-1)/2).
$$
In the rational case, $C_{i_1...i_N}$ is nonzero only if $i_1,...,i_N$ 
are distinct, and equals the sign of the permutation 
$[1,...,N]\to [i_1,...,i_N]$. In this case, the formula 
for quantum determinant is formula (3.11) in \cite{EK3} (with $u$ replaced 
by $u/N$). 

This motivates the following definition. Let $A$ be any algebra 
and $X(u)\in \End(k^N)\o A((u))[[h]]$. In this case, define 
$$
\qdet_R(X(u))=\sum_{i_1,...,i_N}C_{i_1,...,i_N}X_{1i_1}(u-h(N-1)/2)...
X_{Ni_N}(u+h(N-1)/2).
$$
We will need this definition below. 

\subhead 1.4. The normalized R-matrix\endsubhead  

For any $f\in k(\Sigma)[[h]]$, $f=\sum f_nh^n$,  
and $z\in \Sigma(k[[h]])$, we say that $f$ is regular at $z$
if $f_i$ are regular at $z_0$, where $z_0$ is the reduction of 
$z$ mod $h$.  

For any $a\in \Sigma(k[[h]])$ where 
$R$ and $R^{-1}$ are regular, 
let $\tilde E(a)$ be the shifted basic module over $F(R)$: 
as a vector space $\tilde E(a)=k^N$, and 
$\pi_{\tilde E(a)}(T(u))=R(u-a)$.

{\bf Remark.} The crossing symmetry equations (1.3) have the following 
interpretation: the double dual of the basic module $E(0)$ over $F(R)$ is 
isomorphic to the tensor product of the module $E(-Nh)$ and a
1-dimensional module. 

 For any function 
$f\in k(\Sigma)[[h]]$ which is regular at $a$ together with its inverse, 
we can twist 
$\tilde E$ by $f$ and define a new module $\tilde E_f(a)$ by 
$\pi_{\tilde E_f(a)}(T(u))=f(u-a)R(u-a)$. 

\proclaim{Proposition 1.1} There exists a unique
function $f_0=1+O(h)$ such that 
$\pi_{\tilde E_{f_0}(a)}(\qdet_R(T(u)))=1$. 
\endproclaim

\demo{Proof} Let $\pi_{\tilde E(a)}(\qdet_R(T(u)))=\rho(u-a)$, 
where $\rho\in k(\Sigma)[[h]]$. Then $f_0$ is defined by the condition 
$f_0(u-h(N-1)/2-a)...f_0(u+h(N-1)/2-a)=\rho(u-a)$. This equation 
has a unique solution in $k(\Sigma)[[h]]$, which obviously does not depend on 
$a$.  
$\square$\enddemo

We will call $f_0$ the normalizing function for $R$, and denote 
$f_0R$ by $\bR$. We will call $\bR$ the normalized R-matrix. 

{\bf Remark.} The function $f_0$ does not have a simple explicit expression 
even in the rational case. 

\proclaim{Proposition 1.2} The normalized R-matrix satisfies 
the crossing symmetry equations and the unitarity condition:
$$
\gather
(((\bR(u)^{-1})^{t_1})^{-1})^{t_1}=\bR(u+Nh)\\
(((\bR(u)^{-1})^{t_2})^{-1})^{t_2}=
\bR(u+Nh)\\
\bR(u)\bR^{21}(-u)=1\tag 1.6\endgather
$$
\endproclaim

\demo{Proof} It follows from (1.3) that 
$(((\bR(u)^{-1})^{t_2})^{-1})^{t_2}=\phi(u+Nh)\bR(u+Nh)$
for a suitable scalar function $\phi(u)$. 
This implies that $\tilde E_{f_0}(a)^{**}=\tilde E_{f_0\phi}(a-Nh)$. 
However, for any representation $U$ 
of $F(R)$ the quantum determinant acts in $U$ and $U^{**}$ in the same
way, since it is a grouplike element. Thus, since the 
quantum determinant acts trivially on $\tilde E_{f_0}(a)$, 
we get $\phi=1$. This proves the second identity of (1.6). 
The first identity of (1.6) follows immediately, and the third one 
follows from the first two and the fact that $R(u)R^{21}(-u)$ is a scalar. 
$\square$\enddemo

\subhead 1.5. The $\d$-copseudotriangular 
structure on the algebra $F(R)$\endsubhead

In \cite{EK3}, we defined the notion of a copseudotriangular 
structure on a Hopf algebra $A$. 

Recall from \cite{EK3} the construction of a $\d$-copseudotriangular 
structure on $F(R)$, i.e. a bilinear form $B:F(R)\o F(R)\to k(\Sigma)[[h]]$. 
satisfying equations (1.17)-(1.20) of \cite{EK3}. 

Let $\d$ be the derivation of 
$F(R)$ given by $\d T(u)=T'(u)$. 
Let $f\in k(\Sigma)[[h]]$, and $R^f(u)=f(u)R(u)$. 
Then, as we showed in 
\cite{EK3}, the form $B^f$ defined by the formula
$$
\gather
B^f(T^{1,p+q+1}(u_1)\dots T^{p,p+q+1}(u_p),
T^{p+1,p+q+1}(v_1)\dots T^{p+q,p+q+1}(v_q))(y)=\\
\prod_{i=1}^p\prod_{j=q}^1 (R^f)^{i,p+j}(u_i-v_j+y)
\tag 1.7
\endgather
$$
is a $\d$-copseudotriangular structure on $F(R)$. 

We will be especially interested in the case when $f=f_0$, i.e. $R^f=\bR$. 
In this case, we denote $B^f$ by $B^0$.

\proclaim{Proposition 1.3} For any $X\in F(R)$, one has 
$$
B^0(\qdet_R(T(u)),X)(y)=B^0(X,\qdet_R(T(v)))(y)=\e(X).\tag 1.8
$$
\endproclaim

\demo{Proof} Since the quantum determinant is grouplike, 
it is enough to show that 
$$
B^0(\qdet_R(T(u)),T(v))(y)=B^0(T(u),\qdet_R(T(v)))(y)=1.\tag 1.9
$$
By (1.6), $\bR$ satisfies the unitarity condition $\bR(u)\bR^{21}(-u)=1$, 
which implies that $B^0_{13}(u)B^0_{42}(-u)(\Delta(X)\o\Delta(Y))=
\e(X)\e(Y)$. Thus, it suffices to prove only the first equality of (1.9). 
This equality follows by taking the quantum determinant in the identity 
$B^0(T(u),T(v))(y)=\bR(u-v+y)$. 
$\square$\enddemo

\proclaim{Corollary 1.4} The ideal $I\subset F(R)$ generated by the 
relation $\qdet_R(T(u))=1$ belongs to the kernel of $B^0$ on right and left. 
\endproclaim

\subhead 1.5. The Hopf algebra $F_0(R)$\endsubhead

Define the Hopf algebra $F_0(R)$ to be the quotient of 
$F(R)$ by the relation $\qdet_R(T(u))=1$; that is, $F_0(R)=F(R)/I$.  
This is a flat deformation of the algebra $k[SL_N[[t]]]$. 

By Corollary 1.4, the form $B^0$ descends to $F_0(R)$
and defines a copseudotriangular structure on $F_0(R)$. 
Let us denote this structure by $B$. 
It is defined by the formula
$$
\gather
B(T^{1,p+q+1}(u_1)\dots T^{p,p+q+1}(u_p),
T^{p+1,p+q+1}(v_1)\dots T^{p+q,p+q+1}(v_q))(y)=\\
\prod_{i=1}^p\prod_{j=q}^1 \bR^{i,p+j}(u_i-v_j+y)
\tag 1.10
\endgather
$$

\proclaim{Proposition 1.5} The form $B$ is nondegenerate 
(i.e. has trivial kernel) on both sides. 
\endproclaim

\demo{Proof} Let $R_s(u,h)=R(su,sh)$, and $B_s$ be the form on $F_0(R_s)$
constructed as above.  If $B$ is degenerate, 
then $B_s$ is degenerate for all nonzero values of $s$, therefore 
for $s=0$. But for $s=0$ we have $R_s(u)=R_{rat}(u/N)$, so 
$F_0(R)$ is the dual Yangian with opposite product, 
and by Drinfeld's uniqueness theorem 
(\cite{Dr1}; see also Proposition 3.2 in \cite{EK3})
the form $B_0$ is the form on the dual Yangian
with opposite product corresponding to the unique 
pseudotriangular structure $\RR(u)\in Y({\frak sl}_N)
\o Y({\frak sl}_N)$. 

It remains to show that the form $B_0$ is nondegenerate, i.e. that it is 
left-nondegenerate and right-nondegenerate
(the left and the right kernels are zero). 
First of all, 
the right-nondegeneracy follows from 
left-nondegeneracy and the identity 
$B(Y,X)(-u)=B(S(X),Y)(u)$, which is a consequence of the 
identities \linebreak $\RR(u)\RR^{21}(-u)=1$ and $(S\o 1)(\RR)\RR=1$.  

The left-nondegeneracy can be shown, for example, 
as follows: as we pointed out in \cite{EK3}, $B_0$ is the
quantization 
of the standard $\d$-copseudotriangular structure $\beta$ on the Lie bialgebra
$t^{-1}{\frak sl}_N[t^{-1}]$, which is nondegenerate, so $B_0$ is 
left-nondegenerate by Proposition 1.22 in \cite{EK3}. 

Here is another proof of left-nondegeneracy of $B_0$, which does 
not use Proposition 1.22 of \cite{EK3}. 
Recall that the Yangian $Y({\frak sl}_N)$ 
is generated by 
$$t^*(u)=t^*_{-1}u^{-1}+t^*_{-2}u^{-2}+...,$$ 
with defining relations (1.4),(1.5), where instead of $T(u)$ one has 
substituted 
$T^*(u)=1+ht^*(u)=1+T^*_{-1}u^{-1}+T^*_{-2}u^{-2}+...$, 
and an additional relation 
$$\qdet_{R_0}(T^*(u))=1.\tag 1.11
$$ 
The bilinear form $B_0$ defines a Hopf algebra homomorphism 
$\Theta: F_0(R_0)\to Y({\frak sl_N})\o_kk((u))$,
by $\Theta(X)(Y)=B_0(X,Y)$.  
(the tensor product is completed). It is easy to check that this 
homomorphism is defined by the formula $\Theta(T(u))=T^*(v+u)$
(where $T^*(v+u):=\sum T^{*(m)}(v)u^m/m!$). 

From this formula, it is clear that $\Theta$ is  
injective. Indeed, let $G=SL_N$, and $\tilde G$ be the group of 
$G$-valued regular functions on $\Bbb P^1\setminus 0$ which equal to $1$ at 
infinity.  
Quasiclassically, $\Theta$ defines 
the group homomorphism $\theta: \tilde G(k)\to G[[u]](k((v)))$, 
which assigns to every element $g(*)\in \tilde G$ its 
Taylor expansion at a generic point $v$. This homomorphism has 
Zariski dense image (any power series can be 
approximated by polynomials), 
which implies that the corresponding map of function 
algebras is injective. 

Since $\Theta$ is injective, $B_0$ is left-nondegenerate.

The proposition is proved. 
$\square$\enddemo

\subhead 1.6. Square of the antipode \endsubhead 

One has the following well known proposition. 

\proclaim{Proposition 1.6} In $F_0(R)$, one has  
$$
S^2(T(u))=T(u+Nh).\tag 1.12
$$
\endproclaim

\demo{Proof} First of all, observe that the form $B$ is invariant 
under the antipode $S$: $B(S(x),S(y))=B(x,y)$. This fact follows 
from the fact that $B$ satisfies the hexagon axioms.

Thus, 
$$
\gather
B(S(T^{-1}(u)^{0,n+1}),T(v_n)^{n,n+1}...T(v_1)^{1,n+1})(z)=\\
B(T^{-1}(u)^{0,n+1},
S^{-1}(T(v_1)^{1,n+1})...S^{-1}(T(v_n)^{n,n+1}))(z).\tag 1.13
\endgather
$$
The last expression can be easily computed
from the properties of $B$, using that $T(u)T^{-1}(u)=1$. 
The answer is  
$$
\gather
B(S(T^{-1}(u)^{0,n+1}),T(v_n)^{n,n+1}...T(v_1)^{1,n+1})(z)=\\
\theta^2(\bR^{01}(u-v_1+z)...\bR^{0n}(u-v_n+z))=\\
\theta^2(\bR^{01}(u-v_1+z))...\theta^2(\bR^{0n}(u-v_n+z))=\\
\bR^{01}(u+Nh-v_1+z)...\bR^{0n}(u+Nh-v_n+z),\tag 1.14
\endgather
$$
where $\theta(X)=(X^{-1})^{t_1}$.  
Here we have used the crossing symmetry of $\bR$. 

Analogously, we have 
$$
\gather
B(T(u+Nh)^{0,n+1},T(v_n)^{n,n+1}...T(v_1)^{1,n+1})(z)=\\
\bR^{01}(u+Nh-v_1+z)...\bR^{0n}(u+Nh-v_n+z).\tag 1.15
\endgather
$$
By nondegeneracy of $B$, we have (1.12).
$\square$\enddemo

\head 2. The algebra $F_0(R)_\z$ and its finite-dimensional 
representations.\endhead

\subhead 2.1. The definition of $F_0(R)_\z$\endsubhead  

Now let us define 
the factored Hopf algebra $F_0(R)_\z$ with $n$ factors $F_0(R)$,
corresponding to a collection of points $z_1,...,z_n\in \Sigma(k[[h]])$
(cf. \cite{EK3}).

Let $\z=(z_1,...,z_n)$, $z_i\in\Sigma(k[[h]])$ be 
such points that $R(u)^{\pm 1}$ is regular 
at $z_i-z_j$ for $i\ne j$. Denote the space of such $\z$ by $\Sigma_n$. 

\proclaim{Definition} The algebra $F_0(R)_\z$ is 
the h-adic completion of the 
algebra over $k[[h]]$ whose generators are the entries
of the coefficients 
of formal series $T_i(u)^{\pm 1}\in \End(k^N)\o F_0(R)_\z[[u]]$,
$T_i(u)=T^i_0+T^i_1u+...$, $i=1,...,n$
and the defining relations are
$$
\gather
T_i(u)T_i(u)^{-1}=T_i(u)^{-1}T_i(u)=1,\\
R^{12}(u-v+z_i-z_j)T_i^{13}(u)T_j^{23}(v)=T_j^{23}(v)T_i^{13}(u)R^{12}(u-v
+z_i-z_j),\\
\qdet_R(T_i)=1.
\tag 2.1\endgather
$$
\endproclaim
 
This algebra is a Hopf algebra, with the coproduct, 
counit, and antipode defined by
$$
\Delta(T_i(u))=T_i^{12}(u)T_i^{13}(u),\ \e(T_i(u))=
1,\ S(T_i(u))=T_i^{-1}(u)\tag 2.2
$$
(cf. \cite{EK3}).

It is obvious that $F_0(R)_\z=F_0(R)_{\z+\bold a}$, where 
$\bold a=(a,...,a)$, $a\in \Sigma(k[[h]])$. 
In particular, if $n=1$ then $F_0(R)_\z$ does not depend on $\z$, and 
is isomorphic to $F_0(R)$. 

{\bf Remark 1.} We showed in \cite{EK3} that $F_0(R)_\z$ can be identified 
with the tensor product $F_0(R)_{z_1}\o...\o F_0(R)_{z_n}$ (i.e. with 
$F_0(R)^{\o n}$) as a coalgebra, using the map 
$b_1\o...\o b_n\to b_1...b_n$. This factorization will be useful below. 

{\bf Remark 2.} Let $\bR=1-h\bar r+O(h^2)$. Then 
$\bar r\in {\frak {sl}_N}\o {\frak {sl}_N}$. Let  
$SL_N(r)_\z$ be the Poisson group defined in Section 2.8 of 
\cite{EK3}. Then $F_0(R)_\z$ is a quantization of $SL_N(r)_\z$. 

Analogously to $F_0(R)$, the algebra $F_0(R)_\z$ has 
a $\d$-copseudotriangular 
structure \linebreak $B_\z:F_0(R)_\z\o F_0(R)_\z\to k(\Sigma)[[h]]$. 
This structure is defined by the formula
$$
\gather
B_\z(T_{m_1}^{1,p+q+1}(u_1)\dots T_{m_p}^{p,p+q+1}(u_p),
T_{l_1}^{p+1,p+q+1}(v_1)\dots T_{l_q}^{p+q,p+q+1}(v_q))(y)=\\
\prod_{i=1}^p\prod_{j=q}^1 \bR^{i,p+j}(u_i-v_j+z_{m_i}-z_{l_j}+y)
\tag 2.3
\endgather
$$

In the remainder of this Chapter, we discuss modules and comodules 
over $F_0(R)$ and $F_0(R)_\z$. In this paper, 
all modules will be left modules and all 
comodules will be right comodules, so we will 
subsequently omit the words ``left'' and ``right''. 
 
\subhead 2.2. Finite-dimensional representations of $F_0(R)_\z$\endsubhead

By a representation of $F_0(R)_\z$, we mean a  
k-vector space $W$ together with a $k[[h]]$-linear homomorphism
$F_0(R)_\z\to \End W[[h]]$. 

Let us describe a special class of finite dimensional representations
of $F_0(R)_\z$, which we call rational representations. 

\proclaim{Definition} A finite dimensional representation 
$U$ of $F_0(R)_\z$ is said to be rational if there exists 
a matrix function $L(u)\in \End k^N\o \End U\o k(\Sigma)[[h]]$, 
regular and invertible 
at $z_1,...,z_n$, such that the representation is defined by 
the assignment \linebreak $T_i(u)\to L(u+z_i)$.  
\endproclaim

{\bf Remark.} It is obvious that the direct sum and tensor product of rational
representations is a rational representation. 
Further, a subrepresentation and a dual representation of a rational 
representation is rational. 

It is clear that given a rational representation $U$ of $F_0(R)_\z$, 
the corresponding function $L(u)$ is uniquely determined. 
We will denote this function by $L_U(u)$. 

For any rational representation $U$ of $F_0(R)_\z$, 
and any point $a\in \Sigma(k[[h]])$, 
define the shifted representation 
$U(a)$ by the equality $L_{U(a)}(u)=L_U(u-a)$. The shifted representation is 
defined only if $L_U$ is regular and invertible at $z_i-a$. 
 
{\bf Examples} 1. The trivial representation: $U=k$, $L_U=1$. 

2. The basic representation $\tilde E$:  
$U=k^N$, $L_U(u)=\bR(u)$. 

3. The shifted basic representation $\tilde E(a)$:  
$U=k^N$, $L_U(u)=\bR(u-a)$ ($a\in\Sigma(k[[h]])$). 

Since rational representations form a tensor category,
we can construct rational representations of the form 
$\tilde E(a_1)\o...\o \tilde E(a_n)$.

\subhead 2.3. Comodules over $F_0(R)$\endsubhead

In this section we are interested in finite-dimensional
$F_0(R)$-comodules $V$. 

{\bf Examples of  $F_0(R)$-comodules:}

1. The trivial comodule: $V=k$, $\pi^*(1)=1$. 

2. The basic comodule $E$: $V=k^N$, $\pi^*(\v)=T^{21}(0)(1\o \v)$.

3. The shifted basic comodule $E(a)$ ($a\in hk[[h]]$): 
$V=k^N$, 
$$
\pi^*(\v)=T^{21}(a)(1\o \v).\tag 2.4
$$  

 Comodules over $F_0(R)$ form a tensor category. 
 From now on, it will be convenient to us to define 
the tensor product in the way which is opposite to the usual one:
$$
\pi^*_{V\o W}=(m^{op}\o 1\o 1)\circ \sigma_{23}\circ (\pi_V^*\o \pi_W^*).
$$ 

In particular, we have comodules $E(a_1)\o...\o E(a_n)$. 

Now we will define a functor $\phi_u$ 
from the category of  $F_0(R)$-comodules 
to the category of  $F_0(R)$-modules, parametrized by 
points $u\in \Sigma(k[[h]])$ such that $\bR^{\pm 1}(y)$ is regular at $y=u$. 

Let $V$ be a  $F_0(R)$-comodule. 
Define the  $F_0(R)$-module $\phi_u(V)$ to be $V$ as a vector space, 
with the following action of $F_0(R)$: 
$$
X\v=B_{12}(X\o \pi^*(\v))(u).\tag 2.5
$$
The fact that this is indeed an action follows from the identity \linebreak
$B(XY,Z)=B(X\o Y,\Delta(Z))$, which is a basic property of $B$.

The functor $\phi_u$ is a tensor functor. This follows from the property 
that $B(X,YZ)=B(\Delta(X),Z\o Y)$. 

Using the functor $\phi_u$, 
for any $z_1,...,z_n$ such that $R(y)$ is regular and invertible 
at $z_1,..,z_n$
 we can construct a functor 
$\phi_{z_1,...,z_n}$ from the category of finite dimensional comodules
over $F_0(R)$ to the category of rational representations 
of $F_0(R)_\z$. This functor is defined by the condition that 
$\phi_{z_1,...,z_n}(U)|_{F_0(R)_{z_i}}=\phi_{z_i}(U)$. 

{\bf Example.} $\phi_{z_1,...,z_n}(E)=\tilde E$. 

In view of this example, the representation $\tilde E$ will further 
be denoted simply by $E$. 

\subhead 2.4. The R-matrix for comodules\endsubhead

For any two finite dimensional  $F_0(R)$-comodules 
$V,W$, define the R-matrix
$$
R_{VW}(u):V\o W \to V\o W((u))[[h]]\tag  2.6
$$
by the formula 
$$
R_{VW}(u)(\v\o \w)=B_{13}(u)(\pi^*(\v)\o \pi^*(\w))\tag 2.7
$$
The function $R_{VW}(u)$ is defined on $\Sigma$ outside 
of singularities of $\bR(u)$. Note that $R_{EE}(u)=\bR(u)$. 

It follows from the properties of
$B$ that this R-matrix has the following properties:

(i) the hexagon relations:
$$
R_{V_1\o V_2,W}(u)=R_{V_2W}^{23}(u)R_{V_1W}^{13}(u),
R_{V,W_1\o W_2}(u)=R_{VW_1}^{12}(u)R_{VW_2}^{13}(u);
$$

(ii) the unitarity condition: $\sigma R_{VW}(u)\sigma R_{WV}(-u)=1$,
where $\sigma$ is the permutation.

(iii) the operator $\beta_{VW}(a):=\sigma R_{VW}^{-1}(a): V\o W\to W\o V$
is an intertwiner $\phi_u(V(a))\o \phi_u(W)\to \phi_u(W)\o \phi_u(V(a))$ 
for $a\in \Sigma(k[[h]])$ (whenever it is defined).

  \head 3. Dimodules over $F_0(R)_\z$\endhead

\vskip .05in

In Kac-Moody theory, one is interested in the representation theory 
of the (centrally extended) double of the Lie bialgebra of functions
on a punctured curve, rather
than in representation theory of this algebra itself. Similarly, here we are 
interested in representation theory of the quantum double of 
$F_0(R)_\z$. Unfortunately, in general we do not have a convenient
explicit description of this double. 
Therefore, we will use the notion of a dimodule over a Hopf algebra
(a module and comodule simultaneously, with a consistency condition), 
which is equivalent to the notion of a module over the quantum double
of this Hopf algebra. This notion was defined by us in \cite{EK2}, 
but we repeat the main definitions for the reader's convenience. 

\subhead 3.1. Dimodules\endsubhead 

\proclaim{Definition} 

 Let $A$ be a Hopf algebra over $k$. A vector space $X$ is said to be equipped
with the structure of an $A$-dimodule if it is endowed
with two linear maps $\pi: A\o X\to X$, $\pi^*: X\to A\o X$,
such that $\pi$ is an  action of $A$ on $X$ as an algebra,
$\pi^*$ is a  coaction of $A$ on $X$ as a coalgebra,
and they agree according to the formula
(cf \cite{Dr1}, p. 816)
$$
\pi^*\circ \pi=(m_3 \o\pi)
\circ \sigma_{13}\sigma_{24}\circ (S^{-1}\o 1^{\o 4})
\circ (\Delta_3\o\pi^*),\tag 3.1 
$$      
where $m_3:=m\circ (m\o 1)$,
and $\Delta_3:=(\Delta\o 1)\circ\Delta$.
\endproclaim

Similarly one defines the notion of a dimodule in the case when $A,X$ are 
topologically free $k[[h]]$-modules.

 An $A$-dimodule
is called trivial if it is trivial both as a module and a comodule. 

There is an obvious notion of tensor product of modules and comodules
over $A$.
Namely, for any two modules (comodules) $V,W$
$$
\pi_{V\o W}=(\pi_V\o\pi_W)\circ\sigma_{23}\circ (\Delta\o 1\o 1);
\pi_{V\o W}^*=(m^{op}\o 1\o 1)\circ \sigma_{23}\circ (\pi_V^*\o\pi_W^*).
\tag 3.2
$$
The tensor product of dimodules
is just the tensor product of the underlying modules and comodules. 
It follows from \cite{Dr1}, p. 816,
and can be checked by a direct computation, that in this way one indeed
obtains a new dimodule.  

Thus, dimodules over $A$ form a tensor category. 

According to the results of Drinfeld \cite{Dr1}, the category 
of dimodules over $A$ has a natural structure of a braided tensor category.
The braiding is defined by the formula
$$
\beta=\sigma\circ R, R=(\pi\o 1)\circ \sigma_{12}\circ (1\o\pi^*).\tag 3.3
$$
Drinfeld proved that (3.3) satisfies the hexagon relations. 

The meaning of the notion of a dimodule is clarified by the following 
construction (\cite{Dr1}). Let $A$ be a finite dimensional Hopf algebra over 
$k$, and $M$ a dimodule over $A$. Then $M$ has a natural structure
of a  module over the quantum double $D(A)=A\o A^{*op}$. Namely,
$A$ already acts in $M$, and the action of $A^{*op}$ is defined by
$\pi(a^*\o \v)=(a^*\o 1)( \pi^*(\v))$, $a^*\in A^*$, $\v\in M$. 
It can be checked that these actions
are compatible. Thus, we have defined a functor $F$
from the category of $A$-dimodules to the category of $D(A)$-modules.

\proclaim{Proposition 3.1} The functor $F$ with trivial tensor structure
is an equivalence from the  braided tensor category of $A$-dimodules to
the braided tensor category of  $D(A)$-modules.
\endproclaim

\demo{Proof} Follows from \cite{Dr1}. 
$\square$\enddemo

Now let $B\subset A$ be a Hopf subalgebra, 
and $I\subset A$ be a two-sided Hopf ideal. 

\proclaim{Definition} 
A vector space ($k[[h]]$-module) $X$ is said to be equipped
with the structure of an $(B,A,I)$-dimodule if it is endowed
with two linear maps $\pi: B\o X\to X$, $\pi^*: X\to A/I\o X$,
such that $\pi$ is an action of $B$ on $X$ as an algebra,
$\pi^*$ is a  coaction of $A/I$ on $X$ as a coalgebra,
and they agree according to formula (3.1).
\endproclaim

If $I=0$, we will refer to $(B,A,I)$-dimodules as $(B,A)$-dimodules.
 
Analogously to Proposition 3.1, one can show that if $A$ is finite-dimensional,
then the tensor category of $(B,A,I)$-dimodules is equivalent to the 
tensor category of modules over the Hopf subalgebra $B\o (A/I)^{*op}
\subset D(A)$. 

Now let us define the induction functor $\Ind_{(B,A,I)}^{(A,A,I)}$ 
from the category of 
$(B,A,I)$-dimodules
to the category of $(A,A,I)$-dimodules. 
Let $V$ be an $(B,A,I)$-dimodule, with a coaction
$p^*: V\to A/I\o V$.

\proclaim{Proposition 3.2} There exists a unique 
$(A,A,I)$-dimodule $\Ind_{(B,A,I)}^{(A,A,I)} V$ equal to $A\o_BV$ 
as a vector space, such that 
$$
\pi(a\o b\o \v)=ab\o \v, \pi^*(1\o \v)=\sigma_{12}\circ (1\o p^*(\v)),\
a,b\in A,\ \v\in V,\tag 3.4
$$
where $\sigma_{12}$ denotes the permutation. 
\endproclaim

\demo{Proof} The uniqueness is obvious. The existence follows from
the fact that $\Ind_{(B,A,I)}^{(A,A,I)} V$ is 
isomorphic to $\Ind_B^A V$ as an $A$-module. 
The coaction of $A/I$ on $V$ can be
computed using (3.4) and (3.1).
$\square$\enddemo
\vskip .05in

\subhead 3.2. Finite-dimensional dimodules\endsubhead

For $i=1,\dots ,n$, set 
$\hz_i:=(z_1,\dots ,\hat z_i,\dots ,z_n)\in \Sigma_{n-1}$. Then 
 $F_0(R)_{\hat\z_i}$ is the subalgebra in $F_0(R)_\z$ 
 generated by the entries of  $T_j(u)$, $j\ne i$. 

Let $V$ be a finite dimensional comodule over $F_0(R)$.
Consider the functor $\phi_{z_1-z_i,...,z_n-z_i}$ (where $z_i-z_i$ is omitted).
Applying this functor to $V$, we will get a rational representation 
of $F_0(R)_{\hz_i}$. Denote this rational representation by $V(z_i)$. 

The space $V(z_i)$ is equipped with a structure of an $F_0(R)_\z$-comodule,
induced by the embedding $F_0(R)=F_0(R)_{z_i}\to F_0(R)_\z$. Thus, 
$V(z_i)$ is an $F_0(R)_{\hz_i}$-module and $F_0(R)_\z$-comodule, and these
structures are 
connected 
by the formula $a\v=B_{\z,12}(0)(a,\pi^*(\v))$, $a\in F_0(R)_{\hz_i}$. 

{\bf Remark.} Although $B$ has a pole at $0$, the function $B_\z(a,b)(u)$
is regular at $u=0$ (for generic $z_l$) if $a\in F_0(R)_{z_i}$, 
$b\in F_0(R)_{z_j}$, $i\ne j$. 

\proclaim{Proposition 3.3} The space $V(z_i)$, equipped
with the coaction of $F_0(R)_\z$ and an action of $F_0(R)_{\hz_i}$ 
defined above, is an $(F_0(R)_{\hz_i},F_0(R)_\z)$-dimodule. 
\endproclaim

\demo{Proof} The proof amounts to checking identity (3.1). Clearly,
it is enough to check it on elements $a\o \v$, where $\v\in V(z_i)$, and 
$a$ is among the generators of the algebra $F_0(R)_{\hz_i}$. 
In other words, it is enough to check (3.1) on
the formal series $T_j(u)\o \v$ for all $j\ne i$, $\v\in V$. 

In the following computation, indices above a tensor will denote
the numbers of components of the tensor product where this tensor 
appears.

Let the comodule structure of $V(z_i)$ be given by 
$$
\pi^*(\v)=T_V^{21}(1\o \v),\tag 3.5
$$
 where 
$T_V\in \End V\o F_0(R)$. 

 From (3.5), we get
$$
\pi^*(\pi(T_j^{01}(u) \v^2))=T_V^{21}T_j^{02}(u)\v^2
.\tag 3.6
$$
 Computation of the right hand side of (3.1) on $\v$ yields
$$
\gather
RHS=(m_3\o \pi)(S^{-1}(T_j^{03}(u))T_j^{04}(u)T_j^{01}(u)T_V^{52})\v^5)=\\
(m_3\o \pi)(S^{-1}(T_j^{03}(u))T_j^{04}(u)T_V^{52}T_j^{01}(u)\v^5).
\tag 3.7
\endgather
$$
Using the substitution $T_j^{01}(u)=S^{-1}(T_j^{01}(u)^{-1})$
and $T_V=S^{-1}(T_V^{-1})$, we can simplify
(3.7) as follows:
$$
RHS=(1\o S^{-1}\o 1)
(T_j^{01}(u)T_j^{02}(u)(T_V^{21})^{-1}(T_j^{01})^{-1}(u)\v^2).
\tag 3.8
$$
Thus, Proposition 3.3 is a consequence of the following Lemma.

\vskip .1in
{\bf Lemma.} We have the following version of the Yang-Baxter equation:
$$
T_j^{01}(u)T_j^{02}(u)T_V^{12}=T_V^{12}T^{02}(u)T^{01}(u).\tag 3.9
$$
in $\End k^N\o\End V\o F_0(R)_\z$. 
\vskip .1in 

Indeed, transposing components 1 and 2 in (3.9) and multiplying both 
sides by $(T_V^{21})^{-1}$ from the right, we get  
that the RHS of (3.6) and (3.8) are the same. 

{\it Proof of the Lemma.} We have
$$
\gather
T_j^{01}(u)T_j^{02}(u)T_V^{12}=
B_{\z,43}(0)(T_j^{04}(u)T_j^{02}(u)T_V^{13}T_V^{12})=\\
m_{32}B_{\z,54}(0)(T_V^{14}T_V^{12}T_j^{05}(u)T_j^{03}(u))=
m_{35}B_{\z,24}(0)(T_V^{14}T_V^{15}T_j^{02}(u)T_j^{03}(u)).\tag 3.10
\endgather
$$

On the other hand
$$
\gather
T_V^{12}T^{02}(u)T^{01}(u)=
B_{\z,43}(0)(T_V^{12}T^{02}(u)T^{04}(u)T^{13}_V)=\\
m_{42}B_{\z,35}(0)(T_V^{14}T^{15}_VT^{02}(u)T^{03}(u)).\tag 3.11
\endgather
$$

Expressions (3.10) and (3.11) are equal by property (1.18) in \cite{EK3} of 
the copseudotriangular structure $B_\z$. The Lemma and 
Proposition 3.3 are proved. 
$\square$\enddemo

We will call the dimodules $V(z_i)$ the evaluation dimodules. 

{\bf Example.} Let $V=E$ be the basic comodule over $F_0(R)$, defined by the 
formula $T_V=T$. In this case, $V(z_i)=E(z_i)$ is the shifted basic 
representation of $F_0(R)_{\hz_i}$. 

\subhead 3.3. Central extensions\endsubhead

In this section we will define useful central extensions 
$\hat F_0(R)$, $\hat F_0(R)_\z$ of the Hopf algebras 
$F_0(R)$, $F_0(R)_\z$, which will also be Hopf algebras. 

Define $\hat F_0(R):=F_0(R)\o_kk[c]$, where $c$ is central and 
the tensor product is h-adically complete. Equip 
$\hat F_0(R)$ with the coproduct such that
$c$ is a primitive element, and
$$
\tilde\Delta(a)=e^{\frac{h}{2}(c\o \d-\d\o c)}\Delta(a), a\in F_0(R),\tag 3.12
$$
where $\d$ is the derivation defined above. 
Define
the counit on $\hat F_0(R)$ to be the extension of the counit on 
$F_0(R)$.  It is clear that these coproduct 
and counit satisfy the axioms of a bialgebra.

It is easy to show that the antipode
of $F_0(R)$ extends to the antipode on $\hat F_0(R)$ by 
setting $S(c)=-c$. Thus, $\hat F_0(R)$ is a 
Hopf algebra. 

The Hopf algebra $\hat F_0(R)$ is defined by the generators $T(u)$, $c$
with defining relations (1.4) and $[T(u),c]=0$. 
By (3.12), the coproduct in this Hopf algebra is given by 
$$
\Delta(T(u))=T^{12}(u-hc^{3}/2)T^{13}(u+hc^{2}/2).\tag 3.13
$$

Now introduce the multi-point analog of $\hat F_0(R)$ -- the algebra 
$\hat F_0(R)_\z$. To do this, consider the 
$\d$-copseudotriangular structure $\hat B$ on 
$\hat F_0(R)$ obtained by pulling back the $\d$-copseudotriangular
structure from $F_0(R)$ under the natural morphism of Hopf algebras 
$\hat F_0(R)\to F_0(R)$ ($c\to 0$). Define $\hat F_0(R)_\z$ to be the factored 
Hopf algebra $A_\z$ for $A=\hat F_0(R)$ and $\d$-copseudotriangular 
structure $\hat B$ (see \cite{EK3}).  

It is checked directly from the definition, 
analogously to Proposition 3.25 of \cite{EK3}, that $\hat F_0(R)_\z$ is the 
h-adically complete algebra generated by
entries of $T_i(u)^{\pm 1}$, $i=1,\dots,n$, 
and central elements $c_1,\dots,c_n$, 
with the relations
$$
\gather
\bR^{12}(u-v+z_i-z_j-\frac{h}{2}(c_i^{3}-c_j^{3}))
T^{13}_i(u)T_j^{23}(v)=\\
T^{23}_j(v)T_i^{13}(u)
\bR^{12}(u-v+z_i-z_j+\frac{h}{2}(c_i^{3}-c_j^{3})),\tag 3.14
\endgather
$$
$$
T_i(u)T_i^{-1}(u)=T_i^{-1}(u)T_i(u)=1.\tag 3.15
$$

It is easy to see that there exists a unique Hopf algebra structure on 
$\hat F_0(R)_\z$ such that $c_i$ are primitive elements, and
$$
\Delta(T_i(u))=T_i^{12}(u-hc_i^3/2)T_i^{13}(u+hc_i^2/2),
\e(T_i(u))=1, S(T_i(u))=T_i(u)^{-1}.\tag 3.16
$$

Now fix $K\in k[[h]]$. 
Let $\<c\>$ be the ideal in $\hat F_0(R)$ generated by $c$.
Let $V$ be a finite dimensional $F_0(R)$-comodule. 
Then we can extend $V$ to a $(k[c],
\hat F_0(R),\<c\>)$-dimodule by setting 
$$
c|_{V}=K\Id. \tag 3.17
$$

We denote the obtained dimodule by $V_K$, and call it
an evaluation dimodule with central charge $K$. 

It is easy to generalize this construction to the multi-point situation.
Let $C$ be the ideal in $\hat F_0(R)_\z$ generated by $c_1,\dots ,c_n$.
Let $V(z_i)$ be as above. 
Then we can extend $V(z_i)$ to a $(\hat F_0(R)_{\hat \z_i}\o k[c_i],
\hat F_0(R)_\z,C)$-dimodule by setting 
$$
c_i|_{V(z_i)}=K\Id, c_j|_{V(z_i)}=0, j\ne i.\tag 3.18
$$

We denote the obtained dimodule by $V_K(z_i)$, and call it
an evaluation dimodule with central charge $K$. 

\subhead 3.4. Local dimodules \endsubhead

Let $V$ be a finite dimensional $F_0(R)$-comodule, and 
$V_K$ be the corresponding evaluation dimodule 
with central charge $K$. 
Define an infinite-dimensional $(\hat F_0(R),\hat F_0(R),\<c\>)$-dimodule
$$
\hat V_K:=
\Ind_{(k[c],\hat F_0(R),\<c\>)}^{(\hat F_0(R),\hat F_0(R),\<c\>)}V_K.\tag 3.19
$$

{\bf Remark.} The dimodule $\hV_K$ is
 analogous to highest weight
(Weyl) modules over affine Kac-Moody algebras. 

This construction has a multi-point analogue. 
Let  $V_K(z_i)$ be the evaluation dimodule 
with central charge $K$ defined above. 
Define an infinite-dimensional $(\hat F_0(R)_\z,\hat F_0(R)_\z,C)$-dimodule 
$$
\hat V_K(z_i):=\Ind_{(\hat F_0(R)_{\hat\z_i}\o k[c_i],\hat F_0(R)_\z,C)}^
{(\hat F_0(R)_\z,\hat F_0(R)_\z,C)}V_K(z_i).\tag 3.20
$$
We will call $\hV_K(z_i)$ a local dimodule at the point $z_i$. 

Clearly, the space $\hV_K(z_i)$ is naturally isomorphic to 
$\hat V_K$, via the identification 
$\hat F_0(R)_{z_i}\to \hat F_0(R)$. 
Therefore, $\hV_K(z_i)$ is naturally an $F_0(R)$-comodule. 

Let the comodule structure on $\hat V_K$ be defined by the formula 
$\pi^*(\v)=T_{\hat V_K}^{21}(1\o \v)$, where
$T_{\hat V_K}\in Hom(\hat V_K,\hat V_K\o F_0(R))$.  

The action 
of $\hat F_0(R)_{\hz_i}$ on 
$\hV_K(z_i)$ has a convenient description in terms of 
the element $T_{V_K}$. 

\proclaim{Proposition 3.4} The action 
of $\hat F_0(R)_{\hz_i}$ on $\hV_K(z_i)$ is given by the formula
$$
T_j(u)\to B_{34}(z_j-z_i)(T^{13}(u),\hat T_{V_K}^{24}), 
j\ne i.\tag 3.21 
$$
\endproclaim

\demo{Proof} On the top component $V(z_i)\subset \hV_K(z_i)$,
formula (3.19) follows from the definition of the evaluation dimodule 
$V(z_i)$. To go from the case of the top component to the general case,
it is enough to observe that both sides of (3.19) have the same
commutation relations with $T_i(u)$. 
$\square$\enddemo

Proposition 3.4 implies that $\hV_K(z_i)$ is a locally rational representation
of $F_0(R)_{\hz_i}$, i.e. an (h-adically completed) inductive 
limit of finite dimensional rational representations. 

\subhead 3.5. Global dimodules and 
quantum conformal blocks\endsubhead

Now let $V^1,\dots ,V^n$ be 
finite dimensional $F_0(R)$-comodules. 
Define the $(\hat F_0(R)_\z,\hat F_0(R)_\z,C)$-dimodule 
$$
M_K(V^1,\dots ,V^n,\z):=
\hV^1_K(z_1)\o\dots \o \hV^n_K(z_n).\tag 3.22
$$ 
We will call $M_K(V^1,...,V^n)$ a global dimodule.
Sometimes we will shortly
denote it by $M_K(\z)$.

Let $\tilde F_0(R)_\z=\hat F_0(R)_\z/\<c_i=c_j,i\ne j\>$.
It is clear that $\tilde F_0(R)_\z$ is a Hopf algebra.
We will denote by $c$ the image of $c_i$ in this quotient.  
An important property of $\tilde F_0(R)_\z$, which 
$\hat F_0(R)_\z$ does not have, is that $F_0(R)_\z$ is a subalgebra
of $\tilde F_0(R)_\z$.

It is clear that the dimodule $M_K(V^1,...,V^n,\z)$
descends to a $(\tilde F_0(R)_\z,\tilde F_0(R)_\z,\<c\>)$-dimodule. 

Since $F_0(R)_\z\subset \tilde F_0(R)_\z$ is a subalgebra, 
we have an action of $F_0(R)_\z$ on $M_K(\z)$. Define the space of invariant
functionals 
$$
B_K(V^1,\dots ,V^n,\z)=B_K(\z):=\Hom_{F_0(R)_\z}(M_K(\z),k[[h]]).\tag 3.23
$$

The elements of $B_K(\z)$ are quantum analogues 
of conformal blocks in the Wess-Zumino-Witten model.
Therefore we will call them ``quantum conformal blocks''.

We have a natural evaluation map $\xi:B_K(\z)\to (V^1\o...\o V^n)^*$ 
defined by
$$
\xi(f)(\v_1\o...\o \v_n)=f(\v_1\o...\o \v_n), \v_i\in V^i(z_i)\subset 
\hat V^i_K(z_i)
.\tag 3.24
$$

\proclaim{Proposition 3.5} The map $\xi$ is a linear isomorphism.
\endproclaim

\demo{Proof} The statement follows from Frobenius reciprocity. 
$\square$\enddemo

Thus we have constructed a vector bundle $B_K(\z)$
over the space 
$\Sigma_n$  
and a trivialization $\xi$ of this bundle.

\head 4. The dual Hopf algebra and the double of $F_0(R)$\endhead

\subhead 4.1. The Hopf algebra $F_0(R)^{*op}$\endsubhead

Consider the dual 
Hopf algebra $F_0(R)^{*op}$. It is a topological Hopf algebra, 
equipped with a weak topology. Let us define a Laurent series 
$T^*(v)\in Mat_N(k)\o F_0(R)^{*op}((v))$ (where 
$F_0(R)^{*op}((v)):=F_0(R)^{*op}\o_kk((v))$ is the completed tensor product), 
which topologically generates $F_0(R)^{*op}$. 

Let $\Theta: F_0(R)\to F_0(R)^{*op}((v))$ be defined by $\Theta(X)(Y)=B(X,Y)$. 
This map is a Hopf algebra homomorphism.
By properties of $B$, $\Theta$ is injective, and its image 
(tensored with $k((v))$) is dense. 

Define 
$$
T^*(v)=\Theta(T(0))\in F_0(R)^{*op}((v)).\tag 4.1 
$$ 
Since the image of $\Theta$ is dense, we get that the entries of 
$T^*(v)$ topologically generate $F_0(R)^{*op}$. 

Since $\Theta$ is a Hopf algebra map, the series $T^*(v)$ satisfies 
equations (1.4),(1.5) 
(with $T^*$ instead of $T$) and the equation $\qdet_R(T^*(v))=1$.

{\bf Remark.}  The series $T^*(v)$ is always
 infinite in the negative direction. However, it is convergent 
in the weak topology, i.e. $T_{-m}\to 0$ when $m\to +\infty$.  

For Yangians ($R(u)=R_{rat}(u/N)$), 
the series $T^*(u)$ has the form $1+\sum T^*_{-j}u^{-j}$, but in general
it is infinite in both directions. However, the algebras 
$F_0(R)^{*op}$ should have ``the same size'' for all $R$. This means 
that we should expect that nonnegative coefficients of $T^*(v)$ 
are not independent generators, but in fact express via negative coefficients. 
Indeed, we have 

\proclaim{Proposition 4.1} The entries of the matrices $T^*_{-j}$
for $j\ge 1$ topologically generate $F_0(R)^{*op}$ (i.e. 
$T^*(v)$ can be uniquely reconstructed from its ``principal part'')
\endproclaim

\demo{Proof} We argue as in the proof of Proposition 1.5.
Consider the R-matrix $R_s(u,h)=R(su,sh)$. Suppose that 
negative $T^*_i$ generate a smaller algebra than $F_0(R)^{*op}$. 
Then it is so for $F_0(R_s)$ for all $s\ne 0$. But then the same statement 
holds for $s=0$ by deformation argument. Thus, 
to get a contradiction,  
it is enough to prove Proposition 4.1 in the Yangian case.  
But in the Yangian case, the nonnegative coefficients of $T^*$ vanish 
(except for $T_0^*$, which equals $1$), so the negative ones 
generate the algebra by definition. Thus, the proposition is 
proved.
$\square$\enddemo

Now consider the dual Hopf algebra $\hat F_0(R)^{*op}$ to the Hopf 
algebra $\hat F_0(R)$. Let $\d\in h^{-1}\hat F_0(R)^{*op}$ 
be the primitive element 
which acts according to the rule $\d(c)=-1/h$, $\d(T(u))=0$. 
Then we have in $\hat F_0(R)^{*op}$:
$$
[\d,T_*(u)]=\frac{dT_*(u)}{du}.\tag 4.2
$$
This explains the notation $\d$ for this element. 

\subhead 4.2. The double of $F_0(R)$ and the universal R-matrix\endsubhead

Let $\Cal R$ denote the universal R-matrix of the double 
$D(F_0(R))=F_0(R)\o F_0(R)^{*op}$, and $\tilde\Cal R$ the 
universal R-matrix of $D(\hat F_0(R))=\hat F_0(R)\o \hat F_0(R)^{*op}$
(the tensor products are completed). 
Since have a natural decomposition $\hat F_0(R)=F_0(R)\o k[c]$ as algebras, 
we can consider both $\Cal R,\tilde\Cal R$ as elements of 
$\hat F_0(R)\o \hat F_0(R)^{*op}$. 

\proclaim{Proposition 4.2} One has 
$$
\tilde \Cal R=e^{-h(c\o \d)/2}\Cal R e^{-h(c\o\d)/2}.\tag 4.3
$$  
\endproclaim

\demo{Proof} By Drinfeld's theorem \cite{Dr1} of uniqueness of the R-matrix 
for the double, it is enough to check that (4.3) satisfies the axioms 
of a quasitriangular structure. The hexagon axioms are immediate, 
and the fact that $\tR$ transforms the coproduct to the opposite 
coproduct is checked by a straightforward computation.
$\square$\enddemo

\head 5. Quantum Knizhnik-Zamolodchikov connection\endhead

\subhead 5.1. The quantum Sugawara construction\endsubhead

Consider the $(\hat F_0(R),\hat F_0(R),\<c\>)$-dimodule $\hV_K$
defined by (3.19).
This dimodule is automatically equipped with a coaction of
$F_0(R)$, which is obtained from the coaction 
of $\hat F_0(R)$ by setting $c=0$.
Moreover, the algebra $F_0(R)$ embeds naturally in
$\hat F_0(R)$ as a subalgebra, which means that it acts 
in $\hV_K$. Let us denote these actions and coactions 
by $\bar\pi$ and $\bar\pi^*$ (warning: in general they do not combine into 
a dimodule structure). 

Define the operator
$Q: \hV_K\to \hV_K$, given by 
$$
Qx=\bar\pi((\tau_{Kh/2}S^{-1}\o 1)(\bar\pi^*(x))),\tag 5.1
$$
where $\tau_a(T(u)):=T(u-a)$. 
We call the operator $Q$ the Sugawara operator.

Consider the evaluation comodule $V$ and let $T_V\in\End(V)\o F_0(R)$
be such that $\pi^*(\v)=T_V^{21}(1\o \v)$.
For example, if $V=E(a)$, $a\in hk[[h]]$, then $T_V=T(a)$.
 
Let $\v\in V\subset \hV_K(0)$.
Then from the definition of $Q$ we get
$$
Q\v=m_{12}(\tau_{Kh/2}S^{-1}(T^{21}_V)(1\o \v)),\tag 5.2
$$ 
where $m_{12}$ 
denotes the multiplication of components.

 From now on we set $\k=K+N$. 

Let $\alpha$ be the automorphism of $\hat F_0(R)$ defined by 
$\alpha(T(u))=T(u-\k h)$, $\alpha(c)=c$. Denote
by $\hV_K^\alpha$ the dimodule $\hV_K$ twisted by $\alpha$, i.e.  
it is the same as a vector space, and $T(u)|_{\hV_K^\alpha}=
\alpha(T(u))|_{\hV_K}$. 

\proclaim{Proposition 5.1} 
The linear operator $Q$ is an isomorphism of dimodules:
$\hV_K\to \hV_K^{\alpha}$. In particular, 
$$
\gather
QT(u)Q^{-1}=T(u-\k h),\\
QT^*(u)Q^{-1}=T^*(u-\k h),
\tag 5.3\endgather
$$ 
on $\hV_K$.
 
\endproclaim

\demo{Proof} 
Following Drinfeld, consider the quantum Casimir operator on $\hV_K(0)$:
$U=m_{12}((S^{-1}\o 1)(\tR))$.
By Drinfeld's theorem \cite{Dr2}, it satisfies the following equation:
$$
U^{-1}XU=S^2(X), X\in F_0(R).\tag 5.4
$$ 

On the other hand, using Proposition 4.2 and the definitions 
of $\bar \pi,\bar\pi^*$,  
we can express $U$ in terms of $Q$ as follows:
$$
\gather
U=m_{12}((S^{-1}\o 1)(e^{-h(c\o \d)/2}\R e^{-h(c\o \d)/2}))=\\
m_{12}(e^{h(c\o \d)/2}(S^{-1}\o 1)(\R) e^{h(c\o \d)/2})=
Qe^{Kh\d}.\tag 5.5\endgather
$$ 
Combining (5.4) and (5.5), we get
$$
Q^{-1}XQ=e^{-Kh\d}U^{-1}XUe^{Kh\d}=
e^{Kh\d}S^2(X)e^{-Kh\d}.
$$
Using Proposition 1.6, we get 
the desired result.
$\square$\enddemo
\vskip .05in

\subhead 5.2. The braiding\endsubhead

Let $V^1,...,V^n$ be finite-dimensional $F_0(R)$-comodules. 
Let $i,j\in\{1,...,n\}$, $i\ne j$. Recall the notations and definitions of 
Chapter 3. 

\proclaim{Proposition 5.2} There exists a unique isomorphism of dimodules
$\beta_{ij}: \hV_K^i(z_i)\o \hV_K^j(z_j)\to 
\hV_K^j(z_j)\o \hV_K^i(z_i)$, such that 
$\beta_{ij}|_{V^i\o V^j}=\sigma R_{V^iV^j}^{-1}(z_i-z_j)$.
\endproclaim

\demo{Proof} It follows from our definitions that 
the operator $\sigma R_{V^iV^j}(z_i-z_j)^{-1}: V^i(z_i)\o V^j(z_j)\to
V^j(z_j)\o V^i(z_i)$ is an isomorphism of evaluation dimodules. 
As the dimodules $\hV_K^i(z_i)\o \hV_K^j(z_j), 
\hV_K^i(z_j)\o \hV_K^i(z_i)$ are induced from
$V^i(z_i)\o V^j(z_j),
V^j(z_j)\o V^i(z_i)$, this operator extends uniquely to an isomorphism
of dimodules.
$\square$\enddemo

It is clear 
from Section 2.4 
that $\beta_{ij}$ satisfy the braid relations and the unitarity
condition $\beta_{ij}\beta_{ji}=1$. 
We will call them the braiding maps.

Now we compute $\beta_{ij}$ explicitly. 
Recall that the dimodules $\hV_K^i(z_i)$, $\hV_K^j(z_j)$
have the structure of $F_0(R)$-comodules. 

Let $\RR(u)\in F_0(R)^{*op}\o F_0(R)^{*op}((u))$ be the element 
defined by the condition $B(X,Y)(u)=\<Y\o X,\RR^{-1}(u)\>$. 

\proclaim{Proposition 5.3} $\beta_{ij}(\v\o \w)=\sigma \RR(z_i-z_j)(\v\o \w)$.
\endproclaim

\demo{Proof} On the top component $V^i(z_i)\o V^j(z_j)$, the 
operators $\beta_{ij}$ and $\sigma \RR(z_i-z_j)$
 are the same by the definition of $\beta_{ij}$. Thanks to 
Proposition 5.2, it is now enough to show that 
the map $\v\o \w\to \sigma \RR(z_i-z_j)(\v\o \w)$ 
commutes with $T_i^{ij}(u)$ and $T_j^{ij}(u)$.

According to (3.13), on $\hV^i_K(z_i)\o \hV^j_K(z_j)$ we have 
$$
T_j^{0,ij}(v)=T_j^{0,i}(v-Kh/2)T_j^{0,j}(v).\tag 5.6
$$
Since $\hV^i_K(z_i)$ is a locally rational representation 
of $F(R)_{\hz_i}$ obtained by applying of the functor $\phi_{\hz_i}$ to some 
$F_0(R)$-comodule, we obtain from (5.6)
$$
T_j^{0,ij}(v)=T_i^{*0,i}(v-z_i+z_j-Kh/2)T_j^{0,j}(v).\tag 5.7
$$
Thus, we need to check the relation

$$
\RR^{ij}(u)T^{*0i}(v-u-Kh/2)T^{0j}(v)=
T^{0j}(v)T^{*0i}(v-u+Kh/2)\RR^{ij}(u).\tag 5.8
$$
in $Mat_N(k)\o D(F_0(R))_K\o D(F_0(R))_K$, where $D(F_0(R))_K:=
D(F_0(R))/(c=K)$. 
(the relation that $\sigma \RR(z_i-z_j)$ commutes
with $T_j^{ij}$ is (5.8) for $u=z_i-z_j$).

Recall the quasitriangularity property of $\tR$:
$$
\tR^{12}T^{01}(v-Kh/2)T^{02}(v+Kh/2)=
T^{02}(v-Kh/2)T^{01}(v+Kh/2)\tR^{12}.\tag 5.9
$$
Using Proposition 4.2, we see that 
in terms of $\R$, this property can be written as
$$
\R^{12}T^{01}(v-Kh/2)T^{02}(v)=
T^{02}(v)T^{01}(v+Kh/2)\R^{12}.\tag 5.10
$$
Let us apply the map $\Theta(-u)$ defined in Chapter 4 
to the component 1 of (5.10). 
Using (4.1) and the identity $(\Theta(-u)\o 1)(\R)=\RR(u)$, 
we get exactly (5.8).

The fact that $\sigma 
\RR(z_i-z_j)$ commutes with $T_j^{ij}$ is checked similarly. 
The Proposition is proved. 
$\square$\enddemo
\vskip .05in

\subhead 5.3 The quantum Knizhnik-Zamolodchikov connection\endsubhead

Let $\{\be_i\}$ be the standard basis in $k^n$, i.e. 
$(\be_i)_j=\delta_{ij}$.
Let $\Cal E$ be a vector bundle over the space $\Sigma_n$, and $a\in k[[h]]$.
By an $a$-connection on $E$ we mean any collection of maps
$A_i(\z): \Cal E_\z\to \Cal E_{\z-ah\be_i}$, defined for every $\z\in 
\Sigma_n$,
where $\Cal E_\z$ denotes the fiber of $\Cal E$ over $\z$. 
An $a$-connection $\{A_i\}$ is called flat if
$A_j(\z-ah\be_i)A_i(\z)=A_i(\z-ah\be_j)A_j(\z)$. 

In the previous sections we defined an infinite-dimensional vector bundle
$M_K(\z)=\hV^1_K(z_1)\o...\o \hV^n_K(z_n)$ over the configuration space
$\Sigma_n$. This bundle has a natural trivialization, since $\hV^i_K(z_i)=
V^i\o F_0(R)$ as vector spaces. For any $\z,\z'\in \Sigma_n$, denote by
$I_{\z,\z'}: M_K(\z)\to M_K(\z')$ the identification map defined by
this trivialization. 

Now we define a flat $(-\k)$-connection
in the vector bundle $M_K(\z)$. We do it by prescribing,
for any $i\in \{1,...,n\}$, the holonomy operator 
$A_i(\z): M_K(\z)\to M_K(\z+\k h\be_i)$.  

Set 
$$
A_i(\z)=I_{\z,\z+\k h\be_i}Q_i,\tag 5.11
$$
where $Q_i$ is the operator $Q$ described in section 5.1 which acts
in the $i$-th component of the tensor product. It is obvious
that $\{A_i\}$ define a flat $(-\k)$-connection. 
This $(-\k)$-connection is called the quantum Knizhnik-Zamolodchikov
(KZ) connection.

\proclaim{Proposition 5.4} We have the following identities:
$$
\gather
A_i(\z)T_i^{0,1...n}(u)=T_i^{0,1...n}(u-\k h)
A_i(\z),\\
A_i(\z)T_j(u)=T_j(u)A_i(\z), j\ne i.\tag 5.12\endgather
$$
\endproclaim

\demo{Proof}  
We know that $\hV^i_K(z_i)$ is a locally  rational representation 
of $F_0(R)_{\hz_i}$. Let $L(u)$ be the corresponding matrix function. 
Let $\theta_{ij}=1$ if $i>j$,
$-1$ if $i<j$, and $0$ if $i=j$.
Recall that in $\hV_K^1(z_1)\o...\o \hV_K^n(z_n)$ one has
$$
\Delta_p(T_i^{01}(u))=T_i^{01}(u+\frac{Kh\theta_{1i}}{2})...
T_i^{0p}(u+\frac{Kh\theta_{pi}}{2}).\tag 5.13
$$
Therefore,  
by the definition of $\hV_K^i(z_i)$, we have
in $\hV^1_K(z_1)\o...\o \hV^i_K(z_i)\o \hV^n_K(z_n)$:
$$
\gather
T_i^{0,1...n}(u)=\\ L^{01}(u+z_i-z_1-\frac{Kh}{2})...
L^{0,i-1}(u+z_i-z_{i-1}-\frac{Kh}{2})
T_i^{0i}(u)\times \\
L^{0,i+1}
(u+z_i-z_{i+1}+\frac{Kh}{2})...
L^{0n}(u+z_i-z_n+\frac{Kh}{2}).\tag 5.14\endgather
$$
Using the first relation of (5.3), 
we obtain
$$
\gather
A_i(\z)T_i^{0,1...n}(u)=\\ 
L^{01}(u+z_i-z_1-\frac{Kh}{2})...
L^{0,i-1}(u+z_i-z_{i-1}-\frac{Kh}{2})
T_i^{0i}(u-\k h)
\times \\
L^{0,i+1}
(u+z_i-z_{i+1}+\frac{Kh}{2})...
L^{0n}(u+z_i-z_n+\frac{Kh}{2})A_i(\z)=\\
T_i^{0,1...n}(u-\k h)
A_i(\z).\tag 5.15\endgather
$$
This proves the first relation in (5.12).

Similarly, using the second relation in (5.3), 
one deduces the second relation in (5.12).
$\square$\enddemo

\proclaim{Proposition 5.5} Let $f\in B_K(\z+\k h\be_i)$. Then
the functional $A_i^*(\z)f$ on $M_K(\z)$ defined by
$$
Q_i^*f(x_1\o...\o x_i\o...\o x_n):=f(x_1\o...\o Q_ix_i\o...\o x_n), x_j\in
\hV_K^j(z_j).\tag 5.16
$$
belongs to $B_K(\z)$.
\endproclaim

\demo{Proof} The proposition follows immediately from Proposition 5.2
and the definition of $B_K(\z)$. 
$\square$\enddemo

Proposition 5.5 shows that the $(-\k)$-connection $\{A_i\}$ 
gives rise to a well
defined flat $\k$-connection $\{\tilde A_i\}$
on the vector bundle $B_K(\z)$, defined by the formula 
$\tilde A_i(z)=A_i^*(z-\kappa h\be_i)$. This 
$\k$-connection is called the quantum KZ connection on quantum conformal 
blocks. 

\proclaim{Proposition 5.6} The connection $\{A_i\}$ commutes with 
the braiding maps $\beta_{i,i+1}$. That is,
$\beta_{i,i+1}A_i(\z)=A_{i+1}(\sigma_{i,i+1}\z)\beta_{i,i+1}$,
where $\sigma_{i,i+1}$ is the transposition (i,i+1).
\endproclaim

\demo{Proof} Using the second relation of (5.3), we obtain the identity
$Q_i^{-1}\RR^{ij}(u)Q_i=\RR^{ij}(u-\k h)$. Together
with Proposition 5.3, this implies the proposition.  
$\square$\enddemo

\vskip .05in

\subhead 5.4. Quantum KZ equations\endsubhead

Since the bundle $B_K(\z)$ is equipped with a trivialization 
$B_K(\z)\to (V^1\o...\o V^n)^*$, the 
$\k$-connection $\{\tilde A_i\}$ in $B_K(\z)$ 
defines a $\k$-connection $\{\nabla_i(\z)\}$ in the trivial bundle over
$\Sigma_n$ with fiber $V^1\o...\o V^n$ ($\nabla_i(\z)$ acts from the fiber 
over $\z$ to the fiber over $\z-\k h\be_i$). 
This connection is defined by the formula 
$\nabla_i=(\tilde A_i^*)^{-1}$. In this section
we compute the $\k$-connection $\{\nabla_i\}$ explicitly.

\proclaim{Theorem 5.7}
$$
\gather
\nabla_i(\z)=
R_{V^{i-1}V^i}^{i-1,i}
(z_{i-1}-z_i+\k h)...R_{V^1V^i}^{1i}(z_1-z_i+\k h) \times\\
R^{ni}_{V^nV^i}(z_n-z_i)...R^{i+1i}_{V^{i+1}V^i}(z_{i+1}-z_i).
\tag 5.17
\endgather
$$
\endproclaim

\demo{Proof} We have to compute the expression 
$f(\v_1\o...\o Q_i^{-1}\v_i\o...\o \v_n)$ for an invariant functional
$f\in B_K(\z-\k h\be_i)$. 

Denote by $\R_i$ the element $\R$ 
for the algebra $F_0(R)_{z_i}$. Set 
$\R_i(a):=e^{-a\d^1}\R_i e^{a\d^1}$. 
Set $\R_i(Kh/2)=\sum_j a_j\o b_j$. Then we have $Q_i=\sum_j S^{-1}(a_j)b_j$.

Using the fact that $\tR$ satisfies the hexagon relations, 
we deduce
$$
(\Delta_n\o 1)(\R_i(Kh/2))=
(\prod_{j=1}^{i-1}\R_i^{j,n+1}(Kh))
\R_i^{i,n+1}(Kh/2)\prod_{j=i+1}^n\R_i^{j,n+1}
\tag 5.18
$$
in $\hV^1_K(z_1)\o...\o \hV^n_K(z_n)\o F_0(R)^{*op}$. 

Let $\R_{*i}=(S^{-1}\o 1)(\R_i)=(1\o S)(\R_i)$. Then from (5.18) we get
$$
(\Delta_n\o 1)(\R_{*i}(Kh/2))=
(\prod_{j=n}^{i+1}\R_{*i}^{j,n+1})
\R_{*i}^{i,n+1}(Kh/2)\prod_{j=i-1}^{1}\R_{*i}^{j,n+1}(Kh).
\tag 5.19
$$

We do the rest of the computation for $i=n$.
Using the invariance of $f$ and the identity (5.19),
we get
$$
f(Q_n\R_{*n}^{n-1,n}(Kh)...R_{*n}^{1n}(Kh)(\v_1\o...\o \v_n))=
f(\v_1\o...\o \v_n),\tag 5.20
$$
where $\v_i\in V^i$. From (5.20), we obtain
$$
f(\v_1\o...\o Q_n^{-1}\v_n)=
f(\R_{*n}^{n-1,n}(Kh)...\R_{*n}^{1,n}(Kh)\v_1\o...\o \v_n).\tag 5.21
$$
Observe that since $(S\o 1)(\R)=\R^{-1}$, we have
$\R_{*i}=(S^{-2}\o 1)(\R_i^{-1})=\R_i(hN)^{-1}$. Making this substitution, we
bring (5.20) to the form
$$
f(\v_1\o...\o Q_n^{-1}\v_n)=
f(\R_n^{n-1,n}(\k h )^{-1}...\R_n^{1,n}(\k h )^{-1}\v_1\o...\o \v_n).
\tag 5.22
$$

From the definition of $V^j(z_j)$ and the definition of $R_{VW}$ 
it follows that
$$
\R_n^{jn}|_{V^1(z_1)\o...\o V^n(z_n)}=
R_{V^jV^n}^{jn}(z_j-z_n)^{-1}.
$$   
This shows that 
$$
\nabla_n(\z)=R_{V^{n-1}V^n}^{n-1,n}
(z_{n-1}-z_n+\k h)...
R_{V^1V^n}^{1n}(z_1-z_n+\k h),\tag 5.23
$$
which coincides with (5.17).

To go from the case $i=n$ to the general case, 
it is enough to use the invariance of the quantum KZ connection with respect
to the braiding. Indeed, by Proposition 5.6,
$$
\nabla_i(\z)=
\beta_{ii+1}^{-1}...\beta_{n-1n}^{-1}\nabla_n(\z')\beta_{n-1n}...\beta_{ii+1},
\tag 5.24
$$
where $\z':=(z_1,...,z_{i-1},z_{i+1},...,z_n,z_i)$. 
Using Proposition 5.3, and (5.23), we get
$$
\gather
\nabla_i(\z)=
R^{ii+1}_{V^iV^{i+1}}(z_i-z_{i+1}-\k h)
\sigma_{ii+1}...R^{n-1n}_{V^iV^n}(z_i-z_n-\k h)
\sigma_{n-1n}\\
R_{V^nV^i}^{n-1,n}(z_n-z_i+\k h)...R_{V^{i+1}V^i}^{i+1,n}
(z_{i+1}-z_i+\k h)
R_{V^{i-1}V^i}^{i-1,n}
(z_{i-1}-z_i+\k h)...R_{V^1V^i}^{1n}(z_1-z_i+\k h)\\
\sigma_{n-1n}R^{n-1n}_{V^iV^n}(z_i-z_n)^{-1}...\sigma_{ii+1}
R^{ii+1}_{V^iV^{i+1}}(z_i-z_{i+1})^{-1}.\tag 5.25
\endgather
$$
After cancelations we get (5.17).
The theorem is proved.
$\square$\enddemo

Let $F(\z)$ be a section of the trivial bundle over $\Sigma_n$ with fiber 
$V^1\o...\o V^n$, i.e. a function
on $\Sigma_n$ with values in $V^1\o...\o V^n$. We say that $F$ is flat 
with respect to the $\k$-connection $\nabla=(\nabla_i)$ if
for any $i=1,...,n$ $\nabla_i(\z)F(\z)=F(\z-\k h\be_i)$.
The condition for $F$ to be flat is equivalent to the following 
system of difference equations:
$$
\gather
F(\z-\k h\be_i)=
R_{V^{i-1}V^i}^{i-1,i}
(z_{i-1}-z_i+\k h)...R_{V^1V^i}^{1i}(z_1-z_i+\k h) \times\\
R^{ni}_{V^nV^i}(z_n-z_i)...R^{i+1i}_{V^{i+1}V^i}(z_{i+1}-z_i)F(\z).
\tag 5.26\endgather
$$
This system is called the quantum KZ equations.
More precisely, it can be identified 
with the quantum KZ equations of \cite{FR} 
by reversing the $R$-matrix.

\end
\Refs

\ref\by [Ch] I.V.Cherednik \paper On Irreducible representations 
of elliptic quantum R-algebras\jour Soviet Math. Dokl.; \vol 34\issue 3\yr 
1987\endref

\ref\by [Dr1] V.G.Drinfeld \paper Quantum groups \jour 
Proc. Int. Congr. Math. (Berkeley, 1986)\vol 1\pages 798-820\endref 

\ref\by [Dr2] V.G.Drinfeld \paper 
On almost cocommutative Hopf algebras \jour Len. Math.J.
\vol 1\pages 321-342\yr 1990\endref

\ref\by [EK1] P.Etingof and D. Kazhdan\paper Quantization of Lie bialgebras, I,
q-alg 9506005\jour Selecta math. \vol 2\issue 1\yr 1996\pages 1-41\endref

\ref\by [EK2] P.Etingof and D. Kazhdan\paper Quantization of Lie bialgebras, 
II, (revised version) \yr 1996\jour q-alg 9701038\endref

\ref\by [EK3] P.Etingof and D. Kazhdan\paper Quantization of Lie bialgebras, 
III, (revised version), \yr 1996\jour q-alg 9610030\endref

\ref\by [EF] Enriques, B., and Felder, G.\paper Coinvariants of 
Yangian doubles and quantum Knizhnik-Zamolodchikov equations, 
q-alg 9707012\yr 1997\endref


\ref\by [FR] Frenkel, I.B., and Reshetikhin, N.Yu.\paper Quantum affine
algebras and holonomic difference equations
\jour Comm. Math. Phys.\vol
146\pages 1-60\yr 1992\endref

\ref\by [Ha] Hasegawa, K.\paper Crossing symmetry in elliptic 
solutions of the Yang-Baxter equation and a new L-operator 
for Belavin's solution\jour J.Phys. A: 
Math.Gen.\vol 26\yr 1993\pages 3211-3228
\endref

\ref\by [KS] Kazhdan,D., and Soibelman, Y. \paper Representations of
quantum affine algebras \jour Selecta Math.\vol 1\issue 3\yr 1995\endref

\ref\by [KT] Kuroki,G., and Takebe, T.\paper Twisted Wess-Zumino-Witten 
models on elliptic curves\jour q-alg/9612033\yr 1996\endref

\ref\by [TK] Tsuchiya, A., Kanie, Y.\paper Vertex operators in
conformal field theory on $P^1$ and monodromy representations of braid
group\inbook Adv. Stud. Pure Math.\vol 16\pages 297--372\yr 1988\endref 

\ref \by [TUY] Tsuchiya, A., Ueno, K. and Yamada, Y.
\paper Conformal field theory on universal family of satble curves
with gauge symmetries\inbook  Adv. Stud. Pure Math.\vol
19\pages 459--566 \yr 1992\endref
